\documentclass[fleqn]{mat01}
\usepackage{times,mathtimy,amssymb}
\begin{document}
\setcounter{page}{235}
\firstpage{235}

\font\Bbb=msbm10 scaled\magstephalf
\renewcommand\thedefin{\arabic{defin}}
\newcommand{\sign}{\mbox{\rm sign}}
\newcommand{\Discr}{\mbox{\rm Discr}}
\newcommand{\transrow}[1]{\stackrel{#1}{\longrightarrow}}
\def\d{\hbox{\rm d}}
\def\BMO{\hbox{\rm BMO}}
\def\Lip{\hbox{\rm Lip}}
\newtheorem{defn}{\rm DEFINITION}
\newtheorem{remar}{\it Remark}
\newtheorem{lem}{\it Lemma}
\newtheorem{appl}{\it Application}

\title{Multilinear integral operators and mean oscillation}

\markboth{Lanzhe Liu}{Multilinear integral operators and mean oscillation}

\author{LANZHE LIU}

\address{College of Mathematics and Computer,
        Changsha University of Science and Technology,
        Changsha 410077, People's Republic of China\\
\noindent E-mail: lanzheliu@263.net}

\volume{114}

\mon{August}

\parts{3}

\begin{abstract}
In this paper, the boundedness properties for some multilinear operators
related to certain integral operators from Lebesgue spaces to Orlicz
spaces are obtained. The operators include Calder\'on--Zygmund singular
integral operator, fractional integral operator, Littlewood--Paley
operator and Marcinkiewicz operator.
\end{abstract}


\Date{MS received 24 December 2003; revised 12 June 2004}

\maketitle

\section{Introduction and theorems}

Let $b\in \hbox{BMO}(R^n)$ and $T$ be the Calder\'on--Zygmund singular integral
operator. The commutator $[b, T]$ generated by $b$ and $T$ is defined by
$[b, T]f(x)=b(x)Tf(x)-T(bf)(x)$. By a classical result of Coifman {\it et~al}
[6], we know that the commutator is bounded on
$L^p(R^n)$ for $1<p<\infty$. Chanillo [1] proves a similar result when
$T$ is replaced by the fractional integral operators. In [9], the
boundedness properties for the commutators from Lebesgue spaces to
Orlicz spaces are obtained. As the development of Calder\'on--Zygmund
singular integral operators, fractional integral operators and their
commutators (see [7,10,11,15]), multilinear singular integral operators
have been well-studied. In this paper, we are going to consider some
integral operators and their multilinear operator as follows.

Let $m$ be a positive integer and $A$ be a function on $R^n$. We denote that
\begin{equation*}
 R_{m+1}(A;x,y)=A(x)-\sum_{|\alpha|\le m}\frac{1}{\alpha!}D^\alpha A(y)(x-y)^\alpha.
\end{equation*}

\begin{defn}$\left.\right.$\vspace{.3pc}

{\rm \noindent Let $T\!\!: S\to S'$ be a linear operator and there exists a locally integrable function $K(x,y)$
 on $R^n \times R^n$ such that
\begin{equation*}
 Tf(x)=\int_{R^n}K(x, y)f(y)\d y
\end{equation*}
for every bounded and compactly supported function $f$, where $K$ satisfies, for fixed $\varepsilon>0$ and $\delta\ge 0$,
\begin{equation*}
  |K(x,y)|\le C|x-y|^{-n+\delta}
\end{equation*}
and
\begin{equation*}
  |K(y,x)-K(z,x)|\le C|y-z|^\varepsilon|x-z|^{-n-\varepsilon+\delta},
\end{equation*}
if $2|y-z|\le |x-z|$. The multilinear operator related to the integral operator $T$ is defined by
\begin{equation*}
     T^A(f)(x)=\int\frac{R_{m+1}(A;x,y)}{|x-y|^m}K(x,y)f(y)\d y.
\end{equation*}}
\end{defn}

\begin{defn}$\left.\right.$\vspace{.3pc}

\noindent {\rm Let $F_t(x,y)$ define on $R^n\times R^n\times [0,+\infty)$. Hence we denote that
\begin{equation*}
  F_t(f)(x)=\int_{R^n}F_t(x,y)f(y)\d y
\end{equation*}
for every bounded and compactly supported function $f$ and
\begin{equation*}
 F_t^A(f)(x)=\int_{R^n}\frac{R_{m+1}(A;x,y)}{|x-y|^m}F_t(x,y)f(y)\d y.
\end{equation*}
Let $H$ be the Banach space and $H=\{h\!:\|h\|<\infty\}$. For each fixed $x\in R^n$, we view $F_t(f)(x)$ and $F_t^A(f)(x)$ as a
 mapping from $[0,+\infty)$ to $H$. Then, the multilinear operators related to $F_t$ is defined by
\begin{equation*}
  S^A(f)(x)=\|F_t^A(f)(x)\|,
\end{equation*}
where $F_t$ satisfies, for fixed $\varepsilon>0$ and $\delta\ge 0$,
\begin{equation*}
 \|F_t(x,y)\|\le C|x-y|^{-n+\delta}
\end{equation*}
and
\begin{equation*}
  \|F_t(y,x)-F_t(z,x)\|\le C|y-z|^\varepsilon|x-z|^{-n-\varepsilon+\delta},
\end{equation*}
if $2|y-z|\le |x-z|$. We also define that $S(f)(x)=\|F_t(f)(x)\|$.}
\end{defn}

Note that when $m=0$, $T^A$ and $S^A$ are just the commutators of $T$
and $S$ with $A$ (see [9,12,14,18]). When $m>0$, it is
a non-trivial generalization of the commutators. It is well-known that
multilinear operators are of great interest in harmonic analysis and have
been widely studied by many authors [2--5,7,13]. The main purpose
of this paper is to prove the boundedness properties for the multilinear
operators $T^A$ and $S^A$ from Lebesgue spaces to Orlicz spaces.

Let us introduce some notations. Throughout this paper, $Q$ will denote
a cube of $R^n$ with sides parallel to the axes. For any locally
integrable function $f$, the sharp function of $f$ is defined by
\begin{equation*}
     f^{\#}(x)=\sup_{x\in Q}\frac{1}{|Q|}\int_Q |f(y)-f_Q|\d y,
\end{equation*}
 where, and in what follows, $f_Q=|Q|^{-1}\int_Q f(x)\d x$. It is well-known that (see [8])
\begin{equation*}
  f^{\#}(x)\approx\sup_{x\in Q}\inf_{c\in C}\frac{1}{|Q|}\int_Q|f(y)-c|\d y.
\end{equation*}
 For $1\le r<\infty$ and $0\le \delta<n$, let
\begin{equation*}
 M_{\delta,r}(f)(x)=\sup_{x\in Q}\left(\frac{1}{|Q|^{1-r\delta/n}}\int_Q|f(y)|^r\d y\right)^{1/r}.
\end{equation*}
We say that $f$ belongs to $\hbox{BMO}(R^n)$ if $f^{\#}$ belongs to
$L^\infty(R^n)$ and $\|f\|_{\rm BMO}= \|f\|_{L^\infty}$. More generally, let
$\varphi$ be a non-decreasing positive function and define
$\hbox{BMO}_\varphi(R^n)$ as the space of all functions $f$ such that
\begin{equation*}
  \frac{1}{|Q(x,r)|}\int_{Q(x,r)} |f(y)-f_Q|\d y \le C\varphi(r).
\end{equation*}
 For $\beta>0$, the Lipschitz space $\hbox{Lip}_\beta(R^n)$ is the space of functions $f$ such that
\begin{equation*}
 \|f\|_{{\rm Lip}_\beta}=\sup_{x\neq y}|f(x)-f(y)|/|x-y|^\beta<\infty.
\end{equation*}
 For $f$, $m_f$ denotes the distribution function of $f$, that is $m_f(t)=|\{x\in R^n\!: |f(x)|>t\}|$.

Let $\psi$ be a non-decreasing convex function on $R^+$ with
$\psi(0)=0$. $\psi^{-1}$ denotes the inverse function of $\psi$. The
Orlicz space $L_\psi(R^n)$ is defined by the set of functions $f$ such
that $\int \psi(\lambda|f(x)|)\d x<\infty$ for some $\lambda>0$. The norm
is given by $\|f\|_{L_\psi}=\inf_{\lambda>0}\lambda^{-1}(1+\int\psi(\lambda|f(x)|)\d x)$.

We shall prove the following theorems in \S2.

\begin{theorem}[\!]
Let $0\le\delta<n${\rm ,} $1<p<n/\delta$ and $\varphi${\rm ,} $\psi$ be two
non-decreasing positive functions on $R^+$ with
$\varphi(t)=t^{n/p}\psi^{-1}(t^{-n})$ $($or equivalently
$\psi^{-1}(t)=t^{1/p}\varphi(t^{-1/n}))$. Suppose that $\psi$ is convex{\rm ,}
$\psi(0)=0${\rm ,} $\psi(2t)\le C\psi(t)$. Let $T$ be the same as in
Definition~$1$ such that $T$ is bounded from $L^r(R^n)$ to $L^s(R^n)$ for
any $1<r<n/\delta$ and $1/s=1/r-\delta/n$. Then $T^A$ is bounded from
$L^p(R^n)$ to $L_\psi(R^n)$ if $D^\alpha A\in \hbox{\rm BMO}_\varphi(R^n)$ for all
$\alpha$ with\break $|\alpha|=m$.
\end{theorem}

\begin{theorem}[\!]
Let $0\le\delta<n${\rm ,} $1<p<n/\delta$ and $\varphi${\rm ,} $\psi$ be two
non-decreasing positive functions on $R^+$ with
$\varphi(t)=t^{n/p}\psi^{-1}(t^{-n})$ $($or equivalently
$\psi^{-1}(t)=t^{1/p}\varphi(t^{-1/n}))$. Suppose that $\psi$ is convex{\rm ,}
$\psi(0)=0${\rm ,} $\psi(2t)\le C\psi(t)$. Let $S$ be the same as in
Definition~$2$ such that $S$ is bounded from $L^r(R^n)$ to $L^s(R^n)$ for
any $1<r<n/\delta$ and $1/s=1/r-\delta/n$. Then $S^A$ is bounded from
$L^p(R^n)$ to $L_\psi(R^n)$ if $D^\alpha A\in \BMO_\varphi(R^n)$ for all
$\alpha$ with\break $|\alpha|=m$.
\end{theorem}

\renewcommand\theremar{\!}
\begin{remar}$\left.\right.$

{\rm \begin{enumerate}
\renewcommand\labelenumi{{\rm (\roman{enumi})}}
\leftskip .1pc
\item If $\varphi(t)\equiv 1$ and $\psi(t)=t^p$ for $1<p<\infty$, then $T^A$ and $S^A$ are all bounded on
 $L^p(R^n)$ if $D^\alpha A\in \hbox{BMO}(R^n)$ for all $\alpha$ with $|\alpha|=m$.

\item If $\psi(t)=t^q$ and $\varphi(t)=t^{n(1/p-1/q)}$ for
$1<p<q<\infty$, then, by $\BMO_{t^\beta}=\hbox{\rm Lip}_\beta$ (see Lemma 4 of [9]),
$T^A$ and $S^A$ are all bounded from $L^p(R^n)$ to $L^q(R^n)$ if
$D^\alpha A\in \hbox{\rm Lip}_{n(1/p-1/q)}$ for all $\alpha$ with $|\alpha|=m$.
\end{enumerate}}
\end{remar}

\section{Proofs of theorems}

We begin with the following preliminary lemmas.

\begin{lem}\,$[9]$.\
Let $\varphi$ be a non-decreasing positive function on $R^+$ and $\eta$ be an infinitely
 differentiable function on $R^n$ with compact support such that $\int\eta(x)\d x=1$. Denote that
 $b_t(x)=\int_{R^n}b(x-ty)\eta(y)\d y$. Then $\|b-b_t\|_{\rm BMO}\le C\varphi(t)\|b\|_{{\rm BMO}_\varphi}$.
\end{lem}

\begin{lem}\,$[9]$.\
Let $0<\beta<1$ and $\varphi$ be a non-decreasing positive function on $R^+$ or $\beta=1$.
 Then $\|b_t\|_{{\rm Lip}_\beta}\le Ct^{-\beta}\varphi(t)\|b\|_{{\rm BMO}_\varphi}$.
\end{lem}

\begin{lem}\,$[9]$.\
Suppose $1\le p_2<p<p_1<\infty$, $\rho$ is a non-increasing function on $R^+$, $B$ is a linear
 operator such that $m_{B(f)}(t^{1/p_1}\rho(t))\le Ct^{-1}$ if $\|f\|_{L^{p_1}}\le 1$ and $m_{B(f)}(t^{1/{p_2}}\rho(t))
 \le Ct^{-1}$ if $\|f\|_{L^{p_2}}\le 1$. Then $\int_0^\infty m_{B(f)}(t^{1/p}\rho(t))\d t \le C$ if $\|f\|_{L^p}\le (p/p_1)^{1/p}$.
\end{lem}

\begin{lem}\,$[1]$.\
Suppose that $1\le r<p<n/\beta$ and $1/q=1/p-\beta/n$. Then
$\|M_{\beta, r}(f)\|_{L^q}\le C\|f\|_{L^p}$.
\end{lem}

\begin{lem}\,$[15]$.\
Suppose that $1\le r<\infty$ and $b\in \Lip_\beta$. Then
\begin{equation*}
  \|(b-b_Q)f\chi_{2Q}\|_{L^r}\le C|Q|^{1/r}\|b\|_{{\rm Lip}_\beta}M_{\beta, r}(f)(x).
\end{equation*}
\end{lem}

\begin{lem}\,$[4]$.\
Let $A$ be a function on $R^n$ and $D^\alpha A\in L^q(R^n)$
  for all $\alpha$ with $|\alpha|=m$ and some $q>n$. Then
\begin{equation*}
 |R_m(A;x,y)|\le C|x-y|^m\sum_{|\alpha|=m}\left(\frac{1}{|\tilde Q(x,y)|}\int_{\tilde Q(x,y)}|D^\alpha A(z)|^q \d z\right)^{1/q},
\end{equation*}
 where $\tilde Q$ is the cube centered at $x$ and having side length $5\sqrt{n}|x-y|$.
\end{lem}

To prove the theorems of the paper, we need the following:\vspace{.5pc}

\noindent {\it Key Lemma}.\ \ {\it Let $T$ and $S$ be the same as in Definitions $1$ and $2$. Suppose that $Q=Q(x_0, d)$ is a cube
 with {\rm supp} $f\subset(2Q)^c$ and $x,\tilde x\in Q$.\vspace{.5pc}
\begin{enumerate}
\renewcommand\labelenumi{{\rm (\alph{enumi})}}
\item If $0<\delta<n$ and $D^\alpha A\in \BMO(R^n)$ for all $\alpha$ with $|\alpha|=m${\rm ,} then
\begin{align*}
\hskip -1.4pc& |T^A(f)(x)-T^A(f)(x_0)|\\
\hskip -1.4pc&\quad\le C\sum_{|\alpha|=m}\|D^\alpha A\|_{\rm BMO}(M_{\delta, 1}(f)(\tilde x)+M_{\delta, r}(f)(\tilde x))\quad \mbox{for any}\ \ r>1.
\end{align*}

\item If $0<\beta+\delta<n$ and $D^\alpha A\in \Lip_\beta(R^n)$ for all $\alpha$ with $|\alpha|=m${\rm ,} then
\begin{equation*}
\hskip -1.4pc |T^A(f)(x)-T^A(f)(x_0)|\le C\sum_{|\alpha|=m}\|D^\alpha A\|_{{\rm Lip}_\beta}M_{\beta+\delta,1}(f)(\tilde x).
\end{equation*}

\item If $0<\delta<n$ and $D^\alpha A\in \BMO(R^n)$ for all $\alpha$ with $|\alpha|=m${\rm ,} then
\begin{align*}
\hskip -1.4pc& \|F_t^A(f)(x)-F_t^A(f)(x_0)\|\\
\hskip -1.4pc&\quad\le C\sum_{|\alpha|=m}\|D^\alpha A\|_{\rm BMO}(M_{\delta, 1}(f)(\tilde x)+M_{\delta, r}(f)(\tilde x))\quad \mbox{for any}\ \ r>1.
\end{align*}

\item If $0<\beta+\delta<n$ and $D^\alpha A\in \Lip_\beta(R^n)$ for all $\alpha$ with $|\alpha|=m${\rm ,} then
\begin{equation*}
\hskip -1.4pc \|F_t^A(f)(x)-F_t^A(f)(x_0)\|\le C\sum_{|\alpha|=m}\|D^\alpha A\|_{{\rm Lip}_\beta}M_{\beta+\delta,1}(f)(\tilde x).
\end{equation*}
\end{enumerate}}

\begin{proof}
Let $\tilde A(x)=A(x)-\sum_{|\alpha|=m}\frac{1}{\alpha!}(D^\alpha A)_Qx^\alpha$, then
 $R_{m+1}(A;x,y)=R_{m+1}(\tilde A;x,y)$ and $D^\alpha\tilde A=D^\alpha A-(D^\alpha A)_Q$ for $|\alpha|=m$. Suppose
 supp $f\subset (2Q)^c$ and $x,\tilde x\in Q=Q(x_0,d)$. Note that $|x_0-y|\approx |x-y|$ for $y\in (2Q)^c$. We write
\begin{align*}
\hskip -.5pc &T^A(f)(x)-T^A(f)(x_0)=\int_{R^n}\left[\frac{K(x, y)}{|x-y|^m}-\frac{K(x_0,y)}{|x_0-y|^m}\right]R_m(\tilde A; x,y)f(y)\d y \\[3pt]
\hskip -.5pc &\quad +\int_{R^n}\frac{K(x_0, y)f(y)}{|x_0-y|^m}[R_m(\tilde A; x, y)-R_m(\tilde A; x_0, y)]\d y  \\[3pt]
\hskip -.5pc &\quad -\sum_{|\alpha|=m}\frac{1}{\alpha!}\int_{R^n}\left(\frac{K(x, y)(x-y)^\alpha}{|x-y|^m}-\frac{K(x_0, y)(x_0-y)^\alpha}{|x_0-y|^m}\right)D^\alpha\tilde A(y)f(y)\d y \\[3pt]
\hskip -.5pc &:= I+\hbox{\it II}+\hbox{\it III}.
\end{align*}
\noindent {\rm (a)} By Lemma 6 and the following inequality (see [15]), for $b\in \BMO(R^n)$,
\begin{equation*}
 |b_{Q_1}-b_{Q_2}|\le C\log(|Q_2|/|Q_1|)\|b\|_{\rm BMO}\quad \mbox{for}\ \ Q_1 \subset Q_2,
\end{equation*}
 we know that, for $x\in Q$ and $y\in 2^{k+1}Q\!\setminus\! 2^kQ$ with $k\ge 1$,
\begin{align*}
 |R_m(\tilde A;x,y)|&\le C|x-y|^m\sum_{|\alpha|=m}(\|D^\alpha A\|_{\rm BMO}+|(D^\alpha A)_Q-(D^\alpha A)_{Q(x,y)}|) \\[3pt]
 &\le C k|x-y|^m\sum_{|\alpha|=m}\|D^\alpha A\|_{\rm BMO};
\end{align*}
thus
\begin{align*}
|I|&\le C\int_{R^n\setminus 2Q}\left(\frac{|x-x_0|}{|x_0-y|^{m+n+1-\delta}}+\frac{|x-x_0|^\varepsilon}
 {|x_0-y|^{m+n+\varepsilon-\delta}}\right)\\[3pt]
&\quad\,\times |R_m(\tilde A; x,y)| |f(y)|\d y   \\[3pt]
&\le C\sum_{|\alpha|=m}\|D^\alpha A\|_{\rm BMO}\sum_{k=1}^\infty\int_{2^{k+1}Q\setminus2^kQ}k\left(\frac{|x-x_0|}
 {|x_0-y|^{n+1-\delta}}\right.\\[3pt]
&\qquad\qquad\qquad\qquad\ \quad\qquad\qquad\qquad\,\left.+\frac{|x-x_0|^\varepsilon}{|x_0-y|^{n+\varepsilon-\delta}}\right)|f(y)|\d y   \\[3pt]
&\le C\!\sum_{|\alpha|=m}\!\!\|D^\alpha A\|_{\rm BMO}\sum_{k=1}^\infty k(2^{-k}\!+\!2^{-k\varepsilon})
 \left(\!\frac{1}{|2^{k+1}Q|^{1-\delta/n}}\!\int_{2^{k+1}Q}\!|f(y)|\d y\right)    \\[3pt]
&\le C\sum_{|\alpha|=m}\|D^\alpha A\|_{\rm BMO}M_{\delta,1}(f)(\tilde x).
\end{align*}
For $\hbox{\it II}$, by the formula (see [4])
\begin{equation*}
   R_m(\tilde A; x, y)-R_m(\tilde A; x_0, y)=\sum_{|\eta|<m}\frac{1}{\eta!}R_{m-|\eta|}(D^\eta\tilde A; x, x_0)(x-y)^\eta
\end{equation*}
and Lemma 6, we get
\begin{align*}
 |\hbox{\it II}|&\le C\int_{R^n\setminus 2Q}\frac{|R_m(\tilde A; x,y)-R_m(\tilde A; x_0, y)|}{|x_0-y|^{m+n-\delta}}|f(y)|\d y \\[2pt]
 &\le C\sum_{|\alpha|=m}\|D^\alpha A\|_{\rm BMO}\sum_{k=1}^\infty\int_{2^{k+1}Q\setminus2^kQ}
 \frac{|x-x_0|}{|x_0-y|^{n+1-\delta}}|f(y)|\d y    \\[3pt]
 &\le C\sum_{|\alpha|=m}\|D^\alpha A\|_{\rm BMO}M_{\delta,1}(f)(\tilde x).
\end{align*}
For {\it III}, similar to the estimates of {\it I}, we obtain, for any $r>1$ with $1/r+1/r'=1$,
\begin{align*}
|\hbox{\it III}|&\le C\int_{R^n\setminus 2Q}\left(\frac{|x-x_0|}{|x_0-y|^{n+1-\delta}}+\frac{|x-x_0|^\varepsilon}
 {|x_0-y|^{n+\varepsilon-\delta}}\right)\\[3pt]
&\quad\,\times|D^\alpha A(y)-(D^\alpha A)_Q| |f(y)|\d y   \\[3pt]
&\le
 C\sum_{|\alpha|=m}\sum_{k=1}^\infty(2^{-k}+2^{-k\varepsilon})\left(\frac{1}{|2^{k+1}Q|^{1-r\delta/n}}\int_{2^{k+1}Q}|f(y)|^r\d y\right)^{1/r} \\[3pt]
&\quad\, \times\left(\frac{1}{|2^{k+1}Q|}\int_{2^{k+1}Q}|D^\alpha A(x)-(D^\alpha A)_Q|^{r'}\d x\right)^{1/r'}  \\[3pt]
&\le C\sum_{|\alpha|=m}\|D^\alpha A\|_{\rm BMO}M_{\delta, r}(f)(\tilde x).
\end{align*}
Thus
\begin{align*}
&|T^A(f)(x)-T^A(f)(x_0)|\\[5pt]
&\quad\,\le C\sum_{|\alpha|=m}\|D^\alpha A\|_{\rm BMO}(M_{\delta, 1}(f)(\tilde x)+M_{\delta, r}(f)(\tilde x)).
\end{align*}

\noindent {\rm (b)}\ \ By Lemma 6 and the following inequality, for $b\in \Lip_\beta$,
\begin{equation*}
  |b(x)-b_Q|\le \frac{1}{|Q|}\int_Q\|b\|_{{\rm Lip}_\beta}|x-y|^\beta \d y\le \|b\|_{{\rm Lip}_\beta}(|x-x_0|+d)^\beta,
\end{equation*}
we get
\begin{equation*}
 |R_m(\tilde A;x,y)|\le \sum_{|\alpha|=m}\|D^\alpha A\|_{{\rm Lip}_\beta}(|x-y|+d)^{m+\beta},
\end{equation*}
thus
\begin{align*}
|\hbox{\it I}|&\le C\int_{R^n\setminus 2Q}\left(\frac{|x-x_0|}{|x_0-y|^{m+n+1-\delta}}\right.\\[3pt]
&\qquad\qquad\quad\,\left.+\frac{|x-x_0|^\varepsilon}
 {|x_0-y|^{m+n+\varepsilon-\delta}}\right)|R_m(\tilde A; x,y)| |f(y)|\d y   \\[3pt]
&\le C\sum_{|\alpha|=m}\|D^\alpha A\|_{{\rm Lip}_\beta}\sum_{k=1}^\infty\int_{2^{k+1}Q\setminus2^kQ}\!\left(\frac{|x\!-\!x_0|}
 {|x_0-y|^{n+1-\beta-\delta}}\right.\\[3pt]
&\qquad\qquad\qquad\qquad\qquad\qquad\qquad\quad\,\left.+\frac{|x-x_0|^\varepsilon}{|x_0-y|^{n+\varepsilon-\beta-\delta}}\right)|f(y)|\d y    \\[3pt]
&\le C\!\!\sum_{|\alpha|=m}\!\|D^\alpha A\|_{{\rm Lip}_\beta}\sum_{k=1}^\infty(2^{-k}\!+\!2^{-k\varepsilon})
 \frac{1}{|2^{k\!+\!1}Q|^{1-(\beta+\delta)/n}}\int_{2^{k\!+\!1}Q}\!\!|f(y)|\d y    \\[3pt]
&\le C\sum_{|\alpha|=m}\|D^\alpha A\|_{{\rm Lip}_\beta}M_{\beta+\delta,1}(f)(\tilde x), \\[4pt]
|\hbox{\it II}|&\le C\int_{R^n\setminus 2Q}\frac{|R_m(\tilde A; x,y)-R_m(\tilde A; x_0, y)|}{|x_0-y|^{m+n-\delta}}|f(y)|\d y   \\[3pt]
&\le C\sum_{|\alpha|=m}\|D^\alpha A\|_{{\rm Lip}_\beta}\sum_{k=1}^\infty\int_{2^{k+1}Q\setminus2^kQ}\frac{|x-x_0|}{|x_0-y|^{n+1-\beta-\delta}}|f(y)|\d y \\[3pt]
&\le C\sum_{|\alpha|=m}\|D^\alpha A\|_{{\rm Lip}_\beta}M_{\beta+\delta, 1}(f)(\tilde x),\\[3pt]
|\hbox{\it III}|&\le C\int_{R^n\setminus 2Q}\left(\frac{|x-x_0|}{|x_0-y|^{n+1-\beta-\delta}}+\frac{|x-x_0|^\varepsilon}
 {|x_0-y|^{n+\varepsilon-\beta-\delta}}\right)\\[3pt]
&\quad\,\times|D^\alpha A(y)-(D^\alpha A)_Q| |f(y)|\d y   \\[3pt]
 &\le C\!\!\sum_{|\alpha|=m}\!\!\|D^\alpha A\|_{{\rm Lip}_\beta}\sum_{k=1}^\infty(2^{-k}\!+\!2^{-k\varepsilon})
 \frac{1}{|2^{k+1}Q|^{1-(\beta+\delta)/n}}\int_{2^{k+1}Q}|f(y)|\d y    \\[3pt]
&\le C\sum_{|\alpha|=m}\|D^\alpha A\|_{{\rm Lip}_\beta}M_{\beta+\delta, 1}(f)(\tilde x).
\end{align*}
Thus
\begin{equation*}
 |T^A(f)(x)-T^A(f)(x_0)|\le C\sum_{|\alpha|=m}\|D^\alpha A\|_{{\rm Lip}_\beta}M_{\beta+\delta, 1}(f)(\tilde x).
\end{equation*}

Similar argument as in the proof of (a) and (b) will give the proof of (c) and (d), and so we omit the details.
\end{proof}

Now we are in position to prove our theorems.\vspace{.5pc}

\noindent {\it Proof of Theorem $1$.} \ \ We prove the theorem in several steps. First, we prove
\begin{equation}
 (T^A(f))^{\#}\le C\sum_{|\alpha|=m}\|D^\alpha A\|_{\rm BMO}(M_{\delta, 1}(f)+M_{\delta, r}(f))
\end{equation}
 for any $1<r<n/\delta$. Fix a cube $Q=Q(x_0, d)$ and $\tilde x\in Q $. Let $\tilde A(x)=A(x)-\sum_{|\alpha|=m}\frac{1}{\alpha!}
 (D^\alpha A)_Q x^\alpha$. We write, for $f_1=f\chi_{2Q}$ and $f_2=f\chi_{R^n\setminus 2Q}$,
\begin{align*}
  T^A(f)(x)&=\int_{R^n}\frac{R_{m+1}(A; x, y)}{|x-y|^m}K(x, y)f(y)\d y\\[6pt]
&= \int_{R^n}\frac{R_{m+1}(A; x, y)}{|x-y|^m}K(x, y)f_2(y)\d y  \\[6pt]
 &\quad\,+\int_{R^n}\frac{R_m(\tilde A; x, y)}{|x-y|^m}K(x, y)f_1(y)\d y\\[6pt]
 &\quad\,-\sum_{|\alpha|=m}\frac{1}{\alpha!}\int_{R^n}\frac{K(x, y)(x-y)^\alpha}{|x-y|^m}D^\alpha\tilde A(y)f_1(y)\d y,
 \end{align*}
then
\begin{align*}
 &|T^A(f)(x)-T^A(f_2)(x_0)|   \\[5pt]
 &\quad\le \left|T\left(\frac{R_m(\tilde A;x,\cdot)}{|x-\cdot|^m}f_1\right)(x)\right|+\sum_{|\alpha|=m}\frac{1}{\alpha!}
 \left|T\left(\frac{(x-\cdot)^\alpha}{|x-\cdot|^m}D^\alpha\tilde A f_1\right)(x)\right|\\[5pt]
&\qquad\,+|T^A(f_2)(x)-T^A(f_2)(x_0)| \\[5pt]
 &\quad:=  I_1(x)+I_2(x)+I_3(x),
\end{align*}
thus,
\begin{align*}
 &\frac{1}{|Q|}\int_Q|T^A(f)(x)-T^A(f_2)(x_0)|\d x   \\[5pt]
 &\quad\le \frac{1}{|Q|}\int_Q I_1(x)\d x+\frac{1}{|Q|}\int_QI_2(x)\d x+\frac{1}{|Q|}\int_QI_3(x)\d x  \\[5pt]
 &\quad:= I_1+I_2+I_3.
\end{align*}
Now, for $I_1$, if $x\in Q$ and $y\in 2Q$, using  Lemma 6, we get
\begin{equation*}
 R_m (\tilde A; x,y)\le C|x-y|^m\sum_{|\alpha|=m}\|D^\alpha A\|_{\rm BMO}.
\end{equation*}
Thus, by the $(L^r, L^s)$-boundedness of $T$ for $1/s=1/r-\delta/n$ and Holder's inequality, we obtain
\begin{align*}
 I_1 &\le C\sum_{|\alpha|=m}\|D^\alpha A\|_{\rm BMO}\frac{1}{|Q|}\int_Q|T(f_1)(x)|\d x\\[4pt]
&\le C\sum_{|\alpha|=m}\|D^\alpha A\|_{\rm BMO}\|T(f_1)\|_{L^s}|Q|^{-1/s} \\[4pt]
 &\le C\sum_{|\alpha|=m}\|D^\alpha A\|_{\rm BMO}\|f_1\|_{L^r}|Q|^{-1/s}\le C\sum_{|\alpha|=m}\|D^\alpha A\|_{\rm BMO}M_{\delta, r}(f)(\tilde x).
\end{align*}
For $I_2$, taking $q>1$, $l>1$ such that $1/s=1/q-\delta/n$ and denoting $r=ql$, by the $(L^q, L^s)$-boundedness of $T$,
 we gain
\begin{align*}
I_2 &\le \frac{C}{|Q|}\int_Q |T(\sum_{|\alpha|=m}(D^\alpha A-(D^\alpha A)_Q)f_1)(x)|\d x  \\
&\le C\sum_{|\alpha|=m}\left(\frac{1}{|Q|}\int_Q|T((D^\alpha A-(D^\alpha A)_Q)f_1)(x)|^s\d x\right)^{1/s} \\
&\le C|Q|^{-1/s}\sum_{|\alpha|=m}\|(D^\alpha A-(D^\alpha A)_Q)f_1\|_{L^q}   \\
&\le C\sum_{|\alpha|=m}\left(\frac{1}{|Q|}\int_Q|D^\alpha A(y)-(D^\alpha A)_Q|^{ql'}\d y\right)^{1/ql'}\\
&\quad\,\times\left(\frac{1}{|Q|^{1-r\delta/n}}\int_Q|f(y)|^{ql}\d y\right)^{1/ql} \\
&\le C\sum_{|\alpha|=m}\|D^\alpha A\|_{\rm BMO}M_{\delta, r}(f)(\tilde x).
\end{align*}
For $I_3$, by using Key Lemma, we have
\begin{equation*}
 I_3\le C\sum_{|\alpha|=m}\|D^\alpha A\|_{\rm BMO}(M_{\delta, 1}(f)(\tilde x)+M_{\delta, r}(f)(\tilde x)).
\end{equation*}
We now put these estimates together, and taking the supremum over all $Q$ such that $\tilde x\in Q$, we obtain
\begin{equation*}
  (T^A(f))^{\#}(\tilde x)\le C\sum_{|\alpha|=m}\|D^\alpha A\|_{\rm BMO}(M_{\delta, 1}(f)(\tilde x)+M_{\delta,r}(f)(\tilde x)).
\end{equation*}
Thus, taking $1\le r<p<n/\delta$, $1/q=1/p-\delta/n$ and by Lemma 4, we obtain
\begin{align}
 \|T^A(f)\|_{L^q}&\le C\|(T^A(f))^{\#}\|_{L^q}\nonumber\\
&\le C\sum_{|\alpha|=m}\|D^\alpha A\|_{\rm BMO}(\|M_{\delta, 1}(f)\|_{L^q}
 +\|M_{\delta, r}(f)\|_{L^q})\nonumber\\
 &\le C\sum_{|\alpha|=m}\|D^\alpha A\|_{\rm BMO}\|f\|_{L^p}.
\end{align}
Secondly, we prove that, for $D^\alpha A\in \Lip_\beta(R^n)$ with $|\alpha|=m$,
\begin{equation}
 (T^A(f))^{\#}\le C\sum_{|\alpha|=m}\|D^\alpha A\|_{{\rm Lip}_\beta}(M_{\beta+\delta, r}((f))+M_{\beta+\delta, 1}(f))
\end{equation}
for any $1\le r<n/(\beta+\delta)$. In fact, by Lemma 6, we have, for $x\in Q$ and $y\in 2Q$,
\begin{align*}
&|R_m(\tilde A; x,y)|\le C|x-y|^m\\[6pt]
&\quad \times\sum_{|\alpha|=m}\sup_{z\in 2Q}|D^\alpha A(z)-(D^\alpha A)_Q|
 \le C|x-y|^m|Q|^{\beta/n}\sum_{|\alpha|=m}\|D^\alpha A\|_{{\rm Lip}_\beta}
\end{align*}
and by Lemma 5, we have
\begin{align*}
&\|(D^\alpha A-(D^\alpha A)_{2Q})f\chi_{2Q}\|_{L^r}\\[5pt]
&\le C|Q|^{1/r}\|D^\alpha A\|_{{\rm Lip}_\beta}\left(\frac{1}{|Q|^{1-r\beta/n}}\int_Q|f(y)|^r\d y\right)^{1/r},
\end{align*}
by the $(L^r, L^s)$-boundedness of $T$ for $1/s=1/r-\delta/n$, we obtain
\begin{align*}
&\frac{1}{|Q|}\int_Q|T^A(f)(x)- T^A(f)(x_0)|\d x\\
&\quad\le \frac{1}{|Q|}\int_Q \left|T\left(\frac{R_m(\tilde A;x,\cdot)}{|x-\cdot|^m}f_1\right)(x)\right|\d x \\[3pt]
&\qquad\,+\frac{1}{|Q|}\int_Q\sum_{|\alpha|=m}\frac{1}{\alpha!}\left|T\left(\frac{(x-\cdot)^\alpha}{|x-\cdot|^m}D^\alpha\tilde Af_1\right)(x)\right|\d x\\[3pt]
&\qquad\,+\frac{1}{|Q|}\int_Q|T^A(f_2)(x)-T^A(f_2)(x_0)|\d x\\[5pt]
 &\quad\le \sum_{|\alpha|=m}\|D^\alpha A\|_{{\rm Lip}_\beta}\frac{C}{|Q|^{1/s-\beta/n}}\left(\int_Q|T(f_1)(x)|^s\d x\right)^{1/s}\\[5pt]
&\qquad\,+
 \sum_{|\alpha|=m}\left(\frac{C}{|Q|}\int_Q|T(D^\alpha\tilde A f\chi_{2Q})(x)|^s\d x\right)^{1/s} \\[5pt]
&\qquad\, +\frac{1}{|Q|}\int_Q|T^A(f_2)(x)-T^A(f_2)(x_0)|\d x  \\[5pt]
&\quad\le C\sum_{|\alpha|=m}\|D^\alpha A\|_{{\rm Lip}_\beta}\frac{1}{|Q|^{1/r-(\beta+\delta)/n}}\|f_1\|_{L^r}\\[5pt]
 &\qquad\,+\sum_{|\alpha|=m}\frac{C}{|Q|^{1/s}}\left(\int_{R^n}|(D^\alpha A(x)-(D^\alpha A)_Q)f(x)\chi_{2Q}(x)|^r\d x\right)^{1/r}  \\[5pt]
&\qquad\, +\frac{1}{|Q|}\int_Q|T^A(f_2)(x)-T^A(f_2)(x_0)|\d x  \\[5pt]
&\quad\le C\sum_{|\alpha|=m}\|D^\alpha A\|_{{\rm Lip}_\beta}(M_{\beta+\delta, 1}(f)(\tilde x)+M_{\beta+\delta, r}(f)(\tilde x)).
\end{align*}
Thus, taking $1\le r<p<n/(\beta+\delta)$, $1/q=1/p-(\beta+\delta)/n$ and by Lemma 4, we obtain
\begin{align}
 \|T^A(f)\|_{L^q}&\le C\|(T^A(f))^{\#}\|_{L^q}\nonumber\\
&\le C\sum_{|\alpha|=m}
 \|D^\alpha A\|_{{\rm Lip}_\beta}\left(\|M_{\beta+\delta, r}(f)\|_{L^q}+\|M_{\beta+\delta, 1}(f)\|_{L^q}\right)\nonumber\\
 &\le C\sum_{|\alpha|=m}\|D^\alpha A\|_{{\rm Lip}_\beta}\|f\|_{L^p}.
\end{align}
Now we verify  that $T^A$ satisfies the conditions of Lemma 3. In fact, for any  $1<p_i<n/(\beta+\delta)$ with
 $1/s_i=1/p_i-\delta/n$, $1/q_i=1/p_i-(\beta+\delta)/n\ (i=1,2)$ and $\|f\|_{L^{p_i}}\le 1$, note that
 $T^A(f)(x)=T^{A-A_r}(f)(x)+T^{A_r}(f)(x)$ and $D^\alpha(A_r)=(D^\alpha A)_r$. By (2) and Lemma 1, we obtain
\begin{align*}
 \|T^{A-A_r}(f)\|_{L^{s_i}}&\le C\sum_{|\alpha|=m}\|D^\alpha (A-A_r)\|_{\rm BMO}\|f\|_{L^{p_i}}\\[4pt]
&\le C\sum_{|\alpha|=m}\|D^\alpha A-(D^\alpha A)_r\|_{\rm BMO} \\[4pt]
 &\le C\sum_{|\alpha|=m}\|D^\alpha A\|_{{\rm BMO}_\varphi}\varphi(r).
\end{align*}
and by (4) and Lemma 2, we obtain
\begin{equation*}
 \|T^{A_r}(f)\|_{L^{q_i}}\!\le\! C\!\!\sum_{|\alpha|=m}\!\|(D^\alpha A)_r\|_{{\rm Lip}_\beta}\|f\|_{L^{p_i}}
 \!\le\! Cr^{-\beta}\varphi(r)\!\!\sum_{|\alpha|=m}\!\|D^\alpha A\|_{{\rm BMO}_\varphi}.
\end{equation*}
Thus, for $r=t^{-1/n}$,
\begin{align*}
 m_{T^A(f)}(t^{1/s_i}\varphi(t^{-1/n}))&\le m_{T^{A-A_r}(f)}(t^{1/s_i}\varphi(t^{-1/n})/2)\\[6pt]
&\quad\,+m_{T^{A_r}(f)}(t^{1/s_i}\varphi(t^{-1/n})/2) \\[6pt]
 &\le C\left[\left(\frac{\varphi(r)}{t^{1/s_i}\varphi(r)}\right)^{s_i}+\left(\frac{r^{-\beta}\varphi(r)}{t^{1/s_i}\varphi(r)}\right)^{q_i}\right]
 =Ct^{-1}.
\end{align*}
Taking $1<s_2<p<s_1<n/(\beta+\delta)$ and by Lemma 3, we obtain, for $\|f\|_{L^p}\le (p/s_1)^{1/p}$,
\begin{equation*}
 \int_{R^n}\psi(|T^A(f)(x)|)\d x=\int_0^\infty m_{T^A(f)}(\psi^{-1}(t))\d t\le C,
\end{equation*}
thus, $\|T^A(f)\|_{L_\psi}\le C$. This completes the proof of Theorem 1.\vspace{.6pc}

\noindent {\it Proof of Theorem $2$.}\  \  Let $Q$, $\tilde A(x)$, $f_1$ and $f_2$ be the same as in the proof of Theorem 1, we write
\begin{align*}
 F_t^A(f)(x)&= \int_{R^n}\frac{R_{m+1}(\tilde A;x,y)}{|x-y|^m}F_t(x,y)f(y)\d y\\[4pt]
 &= \int_{R^n}\frac{R_{m+1}(\tilde A;x,y)}{|x-y|^m}F_t(x,y)f(y)\d y \\[4pt]
 &\quad\,+\int_{R^n}\frac{R_m(\tilde A;x,y)}{|x-y|^m}F_t(x,y)f_1(y)\d y\\[4pt]
 &\quad\,-\sum_{|\alpha|=m}\frac{1}{\alpha!}
 \int_{R^n}\frac{F_t(x,y)(x-y)^\alpha}{|x-y|^m}D^\alpha\tilde A(y)f_1(y)\d y,
\end{align*}
then
\begin{align*}
&\frac{1}{|Q|}\int_Q|S^A(f)(x)-S^A(f_2)(x_0)|\d x\\[5pt]
&=\frac{1}{|Q|}\int_Q|\|F_t^A(f)(x)\|-\|F_t^A(f_2)(x_0)\| |\d x  \\[5pt]
 &\le \frac{1}{|Q|}\int_Q\left|\left|F_t\left(\frac{R_m(\tilde A;x,\cdot)}{|x-\cdot|^m}f_1\right)(x)\right|\right|\d x\\[5pt]
&\quad\, +\sum_{|\alpha|=m}\frac{1}{\alpha!}\frac{1}{|Q|}\int_Q\left|\left|F_t\left(\frac{(x-\cdot)^\alpha}{|x-\cdot|^m}D^\alpha\tilde A f_1\right)(x)\right|\right|\d x \\[5pt]
 &\quad\, +\frac{1}{|Q|}\int_Q\|F_t^A(f_2)(x)-F_t^A(f_2)(x_0)\|\d x.
\end{align*}
Using the same argument as in the proof of Theorem 1 will give the proof of Theorem~2. Hence we omit the details.

\section{Applications}

In this section we shall apply Theorems 1 and 2 of the paper to some
particular operators such as the Calder\'on--Zygmund singular integral
operator, fractional integral operator, Littlewood--Paley operator and
Marcinkiewicz operator.

\begin{appl}
Calder\'on--Zygmund singular integral operator\vspace{.3pc}

\noindent{\rm  Let $T$ be the Calder\'on--Zygmund operator (see [5,8,16]). The multilinear operator related to $T$ is defined by
\begin{equation*}
     T^A(f)(x)=\int\frac{R_{m+1}(A;x,y)}{|x-y|^m}K(x,y)f(y)\d y.
\end{equation*}
 Then it is easy to verify that Key Lemma holds for $T^A$ with $\delta=0$, and thus $T$ satisfies the conditions
 in Theorem 1. The conclusion of Theorem 1 holds for $T^A$ with $\delta=0$.}
\end{appl}

\begin{appl}
Fractional integral operator with rough kernel\vspace{.3pc}

\noindent {\rm For $0<\delta<n$, let $T_\delta$ be the fractional integral operator with rough kernel defined by (see [1,7])
\begin{equation*}
 T_\delta f(x)=\int_{R^n}\frac{\Omega(x-y)}{|x-y|^{n-\delta}}f(y)\d y.
\end{equation*}
The multilinear operator related to $T_\delta$ is defined by
\begin{equation*}
 T_\delta^A f(x)=\int_{R^n}\frac{R_{m+1}(A;x,y)}{|x-y|^{m+n-\delta}}\Omega(x-y)f(y)\d y,
\end{equation*}
where $\Omega$ is homogeneous of degree zero on $R^n$, $\int_{S^{n-1}}\Omega(x')\d\sigma(x')=0$
 and $\Omega\in \Lip_\varepsilon(S^{n-1})$ for some $0<\varepsilon\le 1$, that is there exists a constant $M>0$
 such that for any $x,y\in S^{n-1}$, $|\Omega(x)-\Omega(y)|\le M|x-y|^\varepsilon$. Then $T_\delta$
 satisfies the conditions in Theorem 1. Thus, the conclusion of Theorem 1 holds for $T_\delta^A$.}
\end{appl}

\begin{appl}
Littlewood--Paley operator\vspace{.3pc}

\noindent{\rm Let $\varepsilon>0$, $n>\delta\ge 0$ and $\psi$ be a  fixed function which satisfies the following properties:

(1) \ \ $|\psi(x)|\le C(1+|x|)^{-(n+1-\delta)}$,

(2) \ \ $|\psi(x+y)-\psi(x)|\le C|y|^\varepsilon(1+|x|)^{-(n+1+\varepsilon-\delta)}$, when $2|y|<|x|$.

The multilinear Littlewood--Paley operator is defined by
\begin{equation*}
  g_\psi^A(f)(x)=\left(\int_0^\infty |F_t^A(f)(x)|^2\frac{\d t}{t}\right)^{1/2},
\end{equation*}
where
\begin{equation*}
 F_t^A(f)(x)=\int_{R^n}\frac{R_{m+1}(A;x,y)}{|x-y|^m}\psi_t(x-y)f(y)\d y
\end{equation*}
and $\psi_t(x)=t^{-n+\delta}\psi(x/t)$ for $t>0$. We write that $F_t(f)=\psi_t\ast f$. We also define that
\begin{equation*}
  g_\psi(f)(x)=\left(\int_0^\infty|F_t(f)(x)|^2\frac{\d t}{t}\right)^{1/2},
\end{equation*}
which is the Littlewood--Paley $g$ function (see [17]).

Let $H$ be the space $H=\{h\!:\|h\|=\left(\int_0^\infty|h(t)|^2 \d t/t\right)^{1/2}<\infty\}$, then, for each
  fixed $x\in R^n$, $F_t^A(f)(x)$ may be viewed as a mapping from $[0,+\infty)$ to $H$, and it is clear that
\begin{equation*}
    g_\psi(f)(x)=\|F_t(f)(x)\|\quad \mbox{\rm and}\quad g_\psi^A(f)(x)=\|F_t^A(f)(x)\|.
\end{equation*}
 It is only to verify that Key Lemma holds for $g_\psi^A$. In fact, for $D^\alpha A\in \BMO(R^n)$ with $|\alpha|=m$,
 we write, for a cube $Q=Q(x_0,d)$ with supp $f\subset(2Q)^c$, $x, \tilde x \in Q=Q(x_0,d)$,
\begin{align*}
&F_t^A(f)(x)-F_t^A(f)(x_0)\\[9pt]
&=\int_{R^n}\left(\frac{\psi_t(x-y)}{|x-y|^m}-\frac{\psi_t(x_0-y)}{|x_0-y|^m}\right)R_m(\tilde A; x, y)f(y)\d y \\[9pt]
&\quad\,+\int_{R^n}\frac{\psi_t(x_0-y)}{|x_0-y|^m}(R_m(\tilde A;x,y)-R_m(\tilde A; x_0,y))f(y)\d y \\[9pt]
&\quad\,-\sum_{|\alpha|=m}\frac{1}{\alpha!}\int_{R^n}\left(\frac{(x-y)^\alpha \psi_t(x-y)}{|x-y|^m}
 -\frac{(x_0-y)^\alpha \psi_t(x_0-y)}{|x_0-y|^m}\right)\\[9pt]
&\quad\,\times D^\alpha \tilde A(y)f(y)\d y  \\[9pt]
&:= J_1+J_2+J_3.
\end{align*}
By the condition of $\psi$ and Minkowski's inequality, we obtain, for any $r>1$,
\begin{align*}
\|J_1\|&\le C\int_{R^n}\frac{|R_m(\tilde A; x,y)| |f(y)|}{|x_0-y|^m}\left[\int_0^\infty\left(\frac{t|x-x_0|}
 {|x_0-y|(t+|x_0-y|)^{n+1-\delta}}\right.\right.\\[3pt]
&\hskip 10pc\left.\left.+\frac{t|x-x_0|^\varepsilon}{(t+|x_0-y|)^{n+1+\varepsilon-\delta}}\right)^2\frac{\d t}{t}\right]^{1/2}\d y \\[3pt]
&\le C\int_{(2Q)^c}\left(\frac{|x-x_0|}{|x_0-y|^{m+n+1-\delta}}+\frac{|x-x_0|^\varepsilon}
 {|x_0-y|^{m+n+\varepsilon-\delta}}\right)\\[3pt]
&\quad\,\ \times |R_m(\tilde A; x,y)||f(y)|\d y   \\[3pt]
&\le C\sum_{|\alpha|=m}\|D^\alpha A\|_{\rm BMO}M_{\delta,1}(f)(\tilde x), \\[3pt]
\|J_2\| &\le C\sum_{|\alpha|=m}\|D^\alpha A\|_{\rm BMO}\sum_{k=1}^\infty\int_{2^{k+1}\setminus2^kQ}\frac{|x-x_0|}{|x_0-y|^{n+1-\delta}}|f(y)|\d y \\[3pt]
&\le C\|D^\alpha A\|_{\rm BMO}M_{\delta,1}(f)(\tilde x), \\[3pt]
\|J_3\|&\le C\sum_{|\alpha|=m}\sum_{k=1}^\infty\int_{2^{k+1}\setminus2^kQ}\left(\frac{|x-x_0|}{|x_0-y|^{n+1-\delta}}+\frac{|x-x_0|^\varepsilon}
 {|x_0-y|^{n+\varepsilon-\delta}}\right)\\[3pt]
&\quad\,\ \times |D^\alpha\tilde A(y)||f(y)|\d y      \\[3pt]
&\le C\sum_{|\alpha|=m}\|D^\alpha A\|_{\rm BMO}M_{\delta,r}(f)(\tilde x).
\end{align*}
 Similarly, for $D^\alpha A\in \Lip_\beta(R^n)$ with $|\alpha|=m$, we get
\begin{equation*}
 \|F_t^A(f)(x)-F_t^A(f)(x_0)\|\le C\sum_{|\alpha|=m}\|D^\alpha A\|_{{\rm Lip}_\beta}M_{\beta+\delta,1}(f)(\tilde x).
\end{equation*}
From the above estimates, we know that Theorem 2 holds for $g_\psi^A$.}
\end{appl}

\begin{appl}
Marcinkiewicz operator\vspace{.4pc}

\noindent{\rm  Let $0\le \delta<n$, $0<\varepsilon\le 1$ and $\Omega$ be homogeneous of degree zero on $R^n$ and $\int_{S^{n-1}}\Omega(x')\d\sigma(x')=0$.
 Assume that $\Omega\in \Lip_\varepsilon(S^{n-1})$, that is there exists a constant $M>0$ such that for any $x,y\in S^{n-1}$,
 $|\Omega(x)-\Omega(y)|\le M|x-y|^\varepsilon$. The multilinear Marcinkiewicz operator is defined by
\begin{equation*}
 \mu_\Omega^A(f)(x)=\left(\int_0^\infty |F_t^A(f)(x)|^2\frac{\d t}{t^3}\right)^{1/2},
\end{equation*}
    where
\begin{equation*}
  F_t^A(f)(x)=\int_{|x-y|\le t}\frac{\Omega(x-y)}{|x-y|^{n-1-\delta}}\frac{R_{m+1}(A;x,y)}{|x-y|^m}f(y)\d y.
\end{equation*}
We write that
\begin{equation*}
 F_t(f)(x)=\int_{|x-y|\le t}\frac{\Omega(x-y)}{|x-y|^{n-1-\delta}}f(y)\d y.
\end{equation*}
We also define that
\begin{equation*}
  \mu_\Omega(f)(x)=\left(\int_0^\infty|F_t(f)(x)|^2\frac{\d t}{t^3}\right)^{1/2},
\end{equation*}
which is the Marcinkiewicz operator (see [18]).

Let $H$ be the space $H=\{h\!:\|h\|=\left(\int_0^\infty|h(t)|^2 \d t/t^3\right)^{1/2}<\infty\}$. Then, it is clear
 that
\begin{equation*}
 \mu_\Omega(f)(x)=\|F_t(f)(x)\|\quad \mbox{\rm and}\quad \ \mu_\Omega^A(f)(x)=\|F_t^A(f)(x)\|.
\end{equation*}
Now, it is only to verify that Key Lemma holds for $\mu_\Omega^A$. In fact, for $D^\alpha A\in \BMO(R^n)$ with $|\alpha|=m$,
 a cube $Q=Q(x_0,d)$ with supp $f\subset(2Q)^c$, $x, \tilde x \in Q=Q(x_0,d)$ and $r>1$, we have
\begin{align*}
&\|F_t^A(f)(x)-F_t^A(f)(x_0)\| \\[9pt]
&\le \left(\int_0^\infty\left|\int_{|x-y|\le t}\frac{\Omega(x-y)R_m(\tilde A; x,y)}{|x-y|^{m+n-1-\delta}}f(y)\d y\right.\right.\\[9pt]
&\quad\,\left.\left. -\int_{|x_0-y|\le t}\frac{\Omega(x_0-y)R_m(\tilde A; x_0 , y)}{|x_0-y|^{m+n-1-\delta}}f(y)\d y\right|^2\frac{\d t}{t^3}\right)^{1/2}  \\[9pt]
&\quad\, +\sum_{|\alpha|=m}\left(\int_0^\infty\left|\int_{|x-y|\le t}\left(\frac{\Omega(x-y)(x-y)^\alpha}{|x-y|^{m+n-1-\delta}}\right.\right.\right.\\[9pt]
&\quad\left.\left.\left.-\int_{|x_0-y|\le t}\frac{\Omega(x_0-y)(x_0-y)^\alpha}{|x_0-y|^{m+n-1-\delta}}\right)D^\alpha\tilde A(y)f(y)\d y\right|^2\frac{\d t}{t^3}\right)^{1/2}\\[9pt]
&\le \left(\int_0^\infty\left[\int_{|x-y|\le t,\ |x_0-y|>t}\frac{|\Omega(x-y)| |R_m(\tilde A; x,y)|}{|x-y|^{m+n-1-\delta}}|f(y)|\d y\right]^2
 \frac{\d t}{t^3}\right)^{1/2}   \\[9pt]
&\quad\,+\left(\!\int_0^\infty\left[\!\int_{|x-y|>t, \ |x_0-y|\le t}\!\!\frac{|\Omega(x_0-y)| |R_m(\tilde A; x_0,y)|}{|x_0-y|^{m+n-1-\delta}}|f(y)|\d y\right]^2
 \frac{\d t}{t^3}\right)^{1/2}   \\[9pt]
&\quad\,+\left(\int_0^\infty\left[\int_{|x-y|\le t, |x_0-y|\le t}\left|\frac{\Omega(x-y)R_m(\tilde A; x,y)}{|x-y|^{m+n-1-\delta}}\right.\right.\right.\\[9pt]
&\quad\left.\left.\left.-\frac{\Omega(x_0-y)R_m(\tilde A; x_0,y)}{|x_0-y|^{m+n-1-\delta}}\right\|f(y)|\d y\right]^2\frac{\d t}{t^3}\right)^{1/2}
\end{align*}
\begin{align*}
&\quad\, +\sum_{|\alpha|=m}\left(\int_0^\infty\left|\int_{|x-y|\le t}\left(\frac{\Omega(x-y)(x-y)^\alpha}{|x-y|^{m+n-1-\delta}}\right.\right.\right.\\[8pt]
&\quad\left.\left.\left. -\int_{|x_0-y|\le t}\frac{\Omega(x_0-y)(x_0-y)^\alpha}{|x_0-y|^{m+n-1-\delta}}\right)D^\alpha\tilde A(y)f(y)\d y\right|^2\frac{\d t}{t^3}\right)^{1/2}  \\[8pt]
&:= L_1+L_2+L_3+L_4
\end{align*}
and
\begin{align*}
 L_1 &\le C\int_{R^n}\frac{|f(y)| |R_m(\tilde A; x,y)|}{|x-y|^{m+n-1-\delta}}\left(\int_{|x-y|\le t<|x_0-y|}\frac{\d t}{t^3}\right)^{1/2} \d y   \\[8pt]
 &\le C\int_{R^n}\frac{|f(y)| |R_m(\tilde A; x,y)|}{|x-y|^{m+n-1-\delta}}\left(\frac{1}{|x-y|^2}-\frac{1}{|x_0-y|^2}\right)^{1/2} \d y   \\[8pt]
 &\le C\int_{(2Q)^c}\frac{|f(y)| |R_m(\tilde A; x,y)|}{|x-y|^{m+n-1-\delta}}\frac{|x_0-x|^{1/2}}{|x-y|^{3/2}}\d y   \\[8pt]
 &\le C\sum_{|\alpha|=m}\|D^\alpha A\|_{\rm BMO}M_{\delta,1}(f)(\tilde x).
\end{align*}
Similarly, we have $L_2\le C\sum_{|\alpha|=m}\|D^\alpha A\|_{\rm BMO}M_{\delta,1}(f)(\tilde x)$.

For $L_3$, by the following inequality (see [18]):
\begin{equation*}
 \left|\frac{\Omega(x-y)}{|x-y|^{n-1-\delta}}-\frac{\Omega(x_0-y)}{|x_0-y|^{n-1-\delta}} \right|\le C\left(\frac{|x-x_0|}{|x_0-y|^{n-\delta}}+
    \frac{|x-x_0|^\gamma}{|x_0-y|^{n-1-\delta+\gamma}}\right),
\end{equation*}
we gain
\begin{align*}
 L_3 &\le C\sum_{|\alpha|=m}\|D^\alpha A\|_{\rm BMO}\int_{(2Q)^c}\left(\frac{|x-x_0|}{|x_0-y|^{n-\delta}}
 +\frac{|x-x_0|^\gamma}{|x_0-y|^{n-1-\delta+\gamma}}\right)\\[6pt]
&\quad\,\times\left(\int_{|x_0-y|\le t, |x-y|\le t}\frac{\d t}{t^3}\right)^{1/2}|f(y)|\d y \\[6pt]
 &\le C\sum_{|\alpha|=m}\|D^\alpha A\|_{\rm BMO}M_{\delta,1}(f)(\tilde x).
\end{align*}
 For $L_4$, similar to the proof of $L_1$, $L_2$ and $L_3$, we obtain
\begin{align*}
  L_4&\le C\sum_{|\alpha|=m}\sum_{k=1}^\infty\int_{2^{k+1}Q\setminus2^kQ}\left(\frac{|x-x_0|}{|x_0-y|^{n+1-\delta}}
  +\frac{|x-x_0|^{1/2}}{|x_0-y|^{n+1/2-\delta}}\right.\\[6pt]
&\hskip 10pc\left.+\frac{|x-x_0|^\gamma}{|x_0-y|^{n+\gamma-\delta}}\right)|D^\alpha\tilde A(y)| |f(y)|\d y \\[6pt]
 &\le C\sum_{|\alpha|=m}\|D^\alpha A\|_{\rm BMO}M_{\delta,r}(f)(\tilde x).
\end{align*}
Similarly, for $D^\alpha A\in \Lip_\beta(R^n)$ with $|\alpha|=m$, we get
\begin{equation*}
 \|F_t^A(f)(x)-F_t^A(f)(x_0)\|\le C\sum_{|\alpha|=m}\|D^\alpha A\|_{{\rm Lip}_\beta}M_{\beta+\delta,1}(f)(\tilde x).
\end{equation*}
Thus, Theorem 2 holds for $\mu_\Omega^A$.}
\end{appl}

\section*{Acknowledgement}
The author would like to express his deep gratitude to the referee for his valuable comments
 and suggestions.
This work was supported by the NNSF Grant No. 10271071.

\end{document}